\documentclass[accepted]{ifacconf}

\usepackage{graphicx}      
\usepackage{natbib}        
\usepackage{times}
\usepackage{amsmath}
\usepackage{amsfonts}       
\usepackage{graphicx}

\usepackage{acro}

\usepackage{tikz}
\usepackage{pgfplots}
 \usetikzlibrary{
    angles,
 	arrows,
 	arrows.meta,
 	shapes,
 	shapes.symbols,
 	shapes.callouts,
 	shapes.geometric,
 	arrows,
 	matrix,
 	backgrounds,
 	positioning,
 	plotmarks,
 	calc,
 	patterns,
 	matrix,
 	decorations.pathreplacing,
 	decorations.pathmorphing,
 	decorations.text,
 	decorations.shapes,
 	decorations.fractals,
 	decorations.markings,
 	spy,
    quotes
 }
\usepgfplotslibrary{fillbetween}
\usepgfplotslibrary{statistics}
\usepgfplotslibrary{groupplots}

\pgfplotsset{
  /pgfplots/confidence box/.style 2 args={
    legend image code/.code={
        \definecolor{steelblue31119180}{RGB}{31,119,180}
        \draw[steelblue31119180,no markers, fill=steelblue31119180, opacity=0.5]
        plot coordinates {
        (-0.1cm,-0.1cm)
        (-0.1cm,0.2cm)
        (0.5cm,0.2cm)
        (0.5cm,-0.1cm)
        (-0.1cm,-0.1cm)
      }
      node[rectangle]{};
    }
  }
}

\DeclareMathOperator\supp{supp}

\newcommand{\diag}{\operatorname{diag}}

\usepackage{tikz}
\usetikzlibrary{shapes,arrows}
\usepackage{pgfplots}
\pgfplotsset{compat=newest}
\usepackage{amsmath, amsfonts}

\acsetup{patch/maketitle = false}
\usepackage{todonotes}
\makeatletter
\renewcommand{\todo}[2][]{\tikzexternaldisable\@todo[#1]{#2}\tikzexternalenable}
\makeatother

\usepackage{dsfont}
\DeclareAcronym{aic}{
    short = AIC,
    long = Akaike Information Criterion 
}

\DeclareAcronym{aicc}{
    short = AICc,
    long = Akaike Information Criterion correction 
}

\DeclareAcronym{bfgs}{
    short = BFGS,
    long = Broyden-Fletcher-Goldfarb-Shanno algorithm
}

\DeclareAcronym{bic}{
    short = BIC,
    long = Bayesian Information Criterion 
}

\DeclareAcronym{cks}{
    short = CKS,
    long = Compositional Kernel Search 
}

\DeclareAcronym{dl}{
    short = DL,
    long = Deep Learning 
}

\DeclareAcronym{gp}{
    short = GP,
    long = Gaussian Process,
    long-plural-form = Gaussian Processes 
}

\DeclareAcronym{hmc}{
    short = HMC, 
    long = Hamiltonian Monte Carlo
}

\DeclareAcronym{kl}{
    short = KL, 
    long = Kullback-Leibler
}

\DeclareAcronym{ks}{
    short = KS,
    long = Kernel Search 
}

\DeclareAcronym{lfm}{
    short = LFM, 
    long = Latent Force Model 
}

\DeclareAcronym{lodegp}{
    short = LODE-GP, 
    long = Linear Ordinary Differential Equation Gaussian Process,
    long-plural-form = Linear Ordinary Differential Equation Gaussian Processes
}

\DeclareAcronym{lti}{
    short = LTI, 
    long = Linear Time Invariant 
}

\DeclareAcronym{map}{
    short = MAP, 
    long = Maximum A Posteriori 
}

\DeclareAcronym{mc}{
    short = MC, 
    long = Monte Carlo 
}

\DeclareAcronym{mcmc}{
    short = MCMC, 
    long = Markov Chain Monte Carlo 
}

\DeclareAcronym{mll}{
    short = MLL, 
    long = Marginal Log Likelihood 
}

\DeclareAcronym{ml}{
    short = ML, 
    long = Machine Learning
}

\DeclareAcronym{mpc}{
    short = MPC, 
    long = Model Predictive Control 
}

\DeclareAcronym{nn}{
    short = NN, 
    long = Neural Network 
}

\DeclareAcronym{nuts}{
    short = NUTS, 
    long = No U-Turn Sampler 
}

\DeclareAcronym{ode}{
    short = ODE, 
    long = Ordinary Differential Equation 
}

\DeclareAcronym{rmse}{
    short = rmse, 
    long = Root Mean Squared Error 
}

\DeclareAcronym{se}{
    short = SE, 
    long = Squared Exponential 
}

\DeclareAcronym{sgd}{
    short = SGD, 
    long = Stochastic Gradient Descent 
}

\DeclareAcronym{skc}{
    short = SKC, 
    long = Structured Kernel Composition 
}

\DeclareAcronym{snf}{
    short = SNF, 
    long = Smith Normal Form 
}

\DeclareAcronym{rkhs}{
    short = RKHS, 
    long = Reproducing Kernel Hilbert Space 
}
\begin{document}
\begin{frontmatter}

\title{Linear ordinary differential equations constrained Gaussian Processes for solving optimal control problems\thanksref{footnoteinfo}} 
\thanks[footnoteinfo]{Authors are listed in alphabetical order. Jörn Tebbe and Andreas Besginow are supported by the SAIL project which is funded by the Ministry of Culture and Science of the State North Rhine-Westphalia under grant no NW21-059C.© 2025 the authors. This work has been accepted to IFAC for publication under a Creative Commons Licence CC-BY-NC-ND.}

\author[First]{Andreas Besginow} 
\author[First]{Markus Lange-Hegermann} 
\author[First]{Jörn Tebbe}

\address[First]{Institute Industrial IT - inIT, OWL University of Applied Sciences and Arts, Lemgo Germany}

\begin{abstract}                
This paper presents an intrinsic approach for addressing control problems with systems governed by linear ordinary differential equations (ODEs). We use computer algebra to constrain a Gaussian Process on solutions of ODEs. We obtain control functions via conditioning on datapoints. Our approach thereby connects Algebra, Functional Analysis, Machine Learning and Control theory. We discuss the optimality of the control functions generated by the posterior mean of the Gaussian Process.
We present numerical examples which underline the practicability of our approach.
\end{abstract}


\end{frontmatter}

\section{Introduction}

Algebra plays a role in articulating concepts in systems and control.
For instance, the classical Kalman decomposition is essentially a statement rooted in linear algebra or the theory of $\mathbb{R}[\partial]$-modules in the continuous case and $\mathbb{R}[\sigma]$-modules in the discrete case.
Early developments in algebraic methods were heavily influenced by the objectives of algebraic geometry, namely to ``define numerical invariants and continuous invariants of algebraic varieties, which allow one to distinguish among nonisomorphic varieties'' \cite[I.8]{hartshorne2013algebraic}.
Similarly, algebraic system theory seeks to define invariants that characterize a system.

A cornerstone of algebraic system theory is behaviorism \citep{willems1997introduction}, which emphasizes understanding and describing systems based solely on their externally observable behavior rather than their internal structure or states.
In this framework, a system is represented as the set of its admissible input-output trajectories, shifting the focus from system construction to its behavior. Algebraic structures, such as vector spaces or modules, are often used to model these behaviors, facilitating rigorous reasoning about system dynamics.
This approach provides a structure-agnostic way of studying system properties that remain invariant under base changes.

Algebraic system theory frequently draws on tools from commutative and homological algebra \citep{pommaret1999algebraic}.
By constructing correspondences between systems and algebraic structures, established invariants can be transferred from algebra to system theory.
Moreover, computational algebra methods, such as Gröbner bases, have been employed to compute such invariants \citep{oberst1990multidimensional,fliess1995flatness,zerz2000topics}.
These methods have spurred the development of novel algorithms for Gröbner basis computation \citep{chyzak2005effective,lange2022boundary} and differential algebra \citep{lange2013thomas,lange2020thomas}.
Such algorithms can e.g.\ be used to refine systems via the grade filtration \citep{barakat2010purity,quadrat2013grade}.
Algebraic techniques have even found applications in control design, such as the use of cylindrical algebraic decomposition to construct controllers based on case distinctions arising from inequalities \citep{Fotiou01112006}.

Despite these successes in system description, behavioral approaches and their algebraic formulations have not led to substantial advancements in system control.

The limitations of algebraic descriptions stem from several critical gaps.
These descriptions do not account for data, particularly noisy measurements from systems.
Algebraic frameworks are suited to handling equations but not the actual functions defining system behavior.
Optimal control is often excluded from algebraic approaches, as function norms are not part of conventional algebraic formulations.
Similarly, algebraic descriptions struggle to incorporate safety constraints.

In this paper, we propose to bridge these gaps by combining the strengths of algebraic system theory with methods from functional analysis, measure theory, probability theory, and machine learning in the context of linear system theory.

In linear system theory, a system is typically modeled using a commutative ring $R$ of linear operators and an $R$-module $M$, referred to as the system module.
When $M$ is represented via a presentation matrix $A \in R^{m_1 \times m_0}$, forming an exact sequence:
\[
R^{1\times m_1}\xrightarrow{\cdot A}R^{1\times m_0}\xrightarrow{} M\xrightarrow{} 0
\]
Here, the rows of the matrix $A$ encode the system equations corresponding to $M$.
The concrete behavior of the abstract system $M$ can be reconstructed in a function space\footnote{One often requires $\mathcal{F}$ to be an injective module or even a cogenerator.} $\mathcal{F}$ by applying the left exact functor $\operatorname{Hom}_R(-, \mathcal{F})$, yielding:
\[
\mathcal{F}^{m_1\times1}\xleftarrow{A\cdot}\mathcal{F}^{m_0\times1}\xleftarrow{} \operatorname{Hom}_R(M,\mathcal{F})\xleftarrow{} 0
\]
Hence, $\operatorname{Hom}_R(M,\mathcal{F})$ is the kernel of $A\cdot$, i.e.\ it contains the admissible behaviors that satisfy the equations given by the rows of $A$:
\[
\operatorname{Hom}_R(M,\mathcal{F})\cong \left\{F\in\mathcal{F}^{m_0}\mid AF=0\right\}=:\operatorname{sol}_\mathcal{F}(A)\mbox{.}
\]
For example, $R = \mathbb{R}[\partial_1, \ldots, \partial_d]$ (the ring of linear partial differential operators) and $\mathcal{F} = C^\infty(\mathbb{R}^d, \mathbb{R})$ (the space of smooth functions) form a typical pairing.

We propose a data structure for $\mathcal{F}$ that satisfies the following properties:
\begin{enumerate}
    \item It incorporates a \emph{probability} distribution (or measure) on $\mathcal{F}^{m_0}$ such that \label{item_measure}
    \item the probability distribution is \emph{concentrated} on \label{item_concentrated}
    \item a \emph{dense} subset (in a suitable topology) of $\operatorname{sol}_\mathcal{F}(A)$, and \label{item_dense}
    \item it enables \emph{inference} on potentially noisy data. \label{item_inference}
\end{enumerate}

This structure reflects the fact that $\operatorname{sol}_\mathcal{F}(A)$ is characterized by the probability \eqref{item_measure} with support that contains only \eqref{item_concentrated} and all \eqref{item_dense} behaviors in a way that  allows conditioning on data \eqref{item_inference}.
Conceptually, the probability distribution serves as a Bayesian prior over admissible behaviors, which can be updated based on observations.

To illustrate the advantages of this approach, we focus on systems defined by linear ordinary differential equations with constant coefficients. We assign a probability distribution to their solution sets using Gaussian processes, constructed via the Smith normal form of the system matrix $A$. We demonstrate how this framework facilitates the explicit construction of control functions using the "control as inference" paradigm.

More broadly, our proposed data structure for $\mathcal{F}$ could support:
\begin{itemize}
    \item Parameter identification: The probability distribution depends on system parameters, enabling optimization of the likelihood for parameter estimation.
    \item Uncertainty quantification: Variances     in posterior distributions provide insights into observational uncertainty.
    \item Safety guarantees: The probability of unsafe behavior can be minimized to ensure system safety \citep{tebbe2024efficiently}.
\end{itemize}

\section{Optimal Control}
\label{sec:optimalcontrol}
Optimal control describes the problem of finding an optimal function in a normed function space $(\mathcal{F}^{m_0}, \Vert \cdot \Vert)$ and is therefore related to functional analysis \citep{clarke2013functional}. To this end, Bellman \citep{bellman1966dynamic} and Pontryagin \citep{pontryagin2018mathematical} provided a strong theory on necessary optimality conditions. We will focus on a certain class of optimal control problems including differential equation and point based constraints \citep{clarke2013functional}.

We assume an asymptotic controllable linear \ac{ode} based system with presentation matrix $A$, combining equations of the state $x \in \mathbb{R}^{n_x}$ and control input $u \in \mathbb{R}^{n_u}$ in one variable $f$.
We consider the optimal control problem given as the minimization of the norm of the function $F$ over a given time horizon $[t_0, t_T]$ 

\begin{subequations}
\begin{align} 
    &\min_{F \in \mathcal{F}^{m_0}} \Vert F \Vert \label{eq:mpc:obj}\\
    \text{s.t. } &F \in \text{sol}_\mathcal{F}(A), \label{eq:mpc:con:ode}\\
    &F(t_i) = f_i, \label{eq:mpc:con:init}
\end{align}
\end{subequations}
with $i=1,\dots,p$. These problems are usually solved via numerical methods such as shooting or collocation methods \citep{lewis2012optimal}. In this work we will present a method which determines the solution of such a problem using methods of machine learning, in particular Gaussian Processes, and computer algebra, in particular the Smith Normal Form. 

\section{Gaussian Processes}
A \acf{gp} \citep{rasmussen2006gaussian} $g(t) \sim \mathcal{GP}(\mu(t), k(t, t'))$ is a probability distribution over functions with the property that all $g(t_1), \ldots, g(t_n)$ are jointly Gaussian.
Such a \ac{gp} is fully characterized by its mean $\mu(t)$ and covariance function $k(t, t')$.
By conditioning a \ac{gp} on a noisy dataset $\mathcal{D} = \{(t_1, f_1), \ldots ,(t_n, f_n)\}$ we have the posterior \ac{gp} defined as
\begin{equation}\label{eq:gaussian_process_posterior_distribution}
    \begin{aligned}
        \mu^* &= \mu(t^*) + K_*^T(K + \sigma_n^2 I)^{-1}f \\
        k^* &= K_{**} - K_*^T(K + \sigma_n^2 I)^{-1} {K_*}
    \end{aligned}
\end{equation}
with covariance matrices $K = (k(t_i, t_j))_{i,j} \in \mathbb{R}^{n \times n}$, $K_* = (k(t_i, t^*_j))_{i, j} \in \mathbb{R}^{n \times m}$ and $K_{**} = (k(t^*_i, t^*_j))_{i, j} \in \mathbb{R}^{m \times m}$ for predictive positions $t^* \in \mathbb{R}^m$ with noise variance $\sigma_n^2$.
This is the most common way of applying \acp{gp} in control theory for regression analysis on time series data. 

Additionally, \acp{gp} can be parameterized in terms of hyperparameters $\theta$, which include the, potentially heteroscedastic, noise variance $\sigma_n^2$, i.e. $\sigma_n^2(t) \in \mathbb{R}^{m_0}$.
Additional hyperparameters are commonly introduced via the \acp{gp} covariance function, for example the \ac{se} covariance function often includes signal variance $\sigma_f^2$ and smoothness parameter $\ell^2$: 
\begin{equation}\label{eq:SE_kernel}
    k_{\text{SE}}(t, t') = \sigma_f^2\exp\left(-\frac{(t-t')^2}{2\ell^2}\right)
\end{equation}

These hyperparameters are trained by maximizing the \acp{gp} \ac{mll}:
\begin{equation} \label{eq:gp:MLL}
    \log p(f|t) = -\frac{1}{2}f^T\left(K+\sigma_n^2I\right)^{-1}f - \frac{1}{2}\log \left( \det \left(K+\sigma_n^2I\right) \right)
\end{equation}
where $I$ is the identity and constant terms are omitted. We obtain a quadratic type error term combined with a regularization term based on the determinant of the regularized covariance matrix.

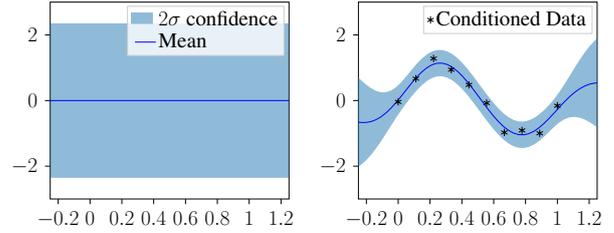
\begin{figure}\label{fig:GP_example}
    \centering
    \scalebox{0.8}{
\begin{tikzpicture}[scale=0.573]

\definecolor{darkgray176}{RGB}{176,176,176}
\definecolor{lightgray204}{RGB}{204,204,204}
\definecolor{steelblue31119180}{RGB}{31,119,180}
\Large
\begin{axis}[
legend cell align={left},
legend style={fill opacity=0.8, draw opacity=1, text opacity=1, draw=lightgray204},
tick align=outside,
tick pos=left,
x grid style={darkgray176},
xmin=-0.25, xmax=1.25,
xtick style={color=black},
y grid style={darkgray176},
ymin=-3, ymax=3,
ytick style={color=black}
]
\path [fill=steelblue31119180, fill opacity=0.5]
(axis cs:-0.25,2.35490489006042)
--(axis cs:-0.25,-2.35490489006042)
--(axis cs:-0.25,-2.35490489006042)
--(axis cs:0.0714285746216774,-2.35490489006042)
--(axis cs:0.157142877578735,-2.35490489006042)
--(axis cs:0.242857158184052,-2.35490489006042)
--(axis cs:0.328571438789368,-2.35490489006042)
--(axis cs:0.414285749197006,-2.35490489006042)
--(axis cs:0.5,-2.35490489006042)
--(axis cs:0.585714280605316,-2.35490489006042)
--(axis cs:0.671428561210632,-2.35490489006042)
--(axis cs:0.757142901420593,-2.35490489006042)
--(axis cs:0.842857122421265,-2.35490489006042)
--(axis cs:0.928571462631226,-2.35490489006042)
--(axis cs:1.0142856836319,-2.35490489006042)
--(axis cs:1.25,-2.35490489006042)
--(axis cs:1.25,2.35490489006042)
--(axis cs:1.25,2.35490489006042)
--(axis cs:1.0142856836319,2.35490489006042)
--(axis cs:0.928571462631226,2.35490489006042)
--(axis cs:0.842857122421265,2.35490489006042)
--(axis cs:0.757142901420593,2.35490489006042)
--(axis cs:0.671428561210632,2.35490489006042)
--(axis cs:0.585714280605316,2.35490489006042)
--(axis cs:0.5,2.35490489006042)
--(axis cs:0.414285749197006,2.35490489006042)
--(axis cs:0.328571438789368,2.35490489006042)
--(axis cs:0.242857158184052,2.35490489006042)
--(axis cs:0.157142877578735,2.35490489006042)
--(axis cs:0.0714285746216774,2.35490489006042)
--(axis cs:-0.0142857134342194,2.35490489006042)
--(axis cs:-0.25,2.35490489006042)
--cycle;

\addlegendentry{$2\sigma$ confidence}
\addlegendimage{
      confidence box={circle,fill,inner sep=1pt}{diamond,fill,inner sep=1pt},
      blue}

\addplot [semithick, blue]
table {%
-0.25 0
-0.0142857134342194 0
0.0714285746216774 0
0.157142877578735 0
0.242857158184052 0
0.328571438789368 0
0.414285749197006 0
0.5 0
0.585714280605316 0
0.671428561210632 0
0.757142901420593 0
0.842857122421265 0
0.928571462631226 0
1.0142856836319 0
1.25 0
};
\addlegendentry{Mean}
\end{axis}

\end{tikzpicture}}
    \scalebox{0.8}{
\begin{tikzpicture}[scale=0.573]

\definecolor{darkgray176}{RGB}{176,176,176}
\definecolor{lightgray204}{RGB}{204,204,204}
\definecolor{steelblue31119180}{RGB}{31,119,180}
\Large
\begin{axis}[
legend cell align={left},
legend style={fill opacity=0.8, draw opacity=1, text opacity=1, draw=lightgray204},
tick align=outside,
tick pos=left,
x grid style={darkgray176},
xmin=-0.25, xmax=1.25,
xtick style={color=black},
y grid style={darkgray176},
ymin=-3, ymax=3,
ytick style={color=black}
]
\addlegendentry{Conditioned Data}
\addplot [semithick, black, mark=asterisk, mark size=3, mark options={solid}, only marks]
table {%
0 -0.0409419648349285
0.111111111938953 0.667891800403595
0.222222223877907 1.28061735630035
0.333333343267441 0.940110445022583
0.444444447755814 0.478916764259338
0.555555582046509 -0.0743785798549652
0.666666626930237 -0.972533524036407
0.777777791023254 -0.910220146179199
0.888888895511627 -0.997552394866943
1 -0.158487573266029
};

\path [fill=steelblue31119180, fill opacity=0.5]
(axis cs:-0.300000011920929,0.934599936008453)
--(axis cs:-0.300000011920929,-2.13949489593506)
--(axis cs:-0.26800000667572,-2.07047033309937)
--(axis cs:-0.236000016331673,-1.96931564807892)
--(axis cs:-0.204000025987625,-1.83424127101898)
--(axis cs:-0.172000020742416,-1.66563546657562)
--(axis cs:-0.140000015497208,-1.46646726131439)
--(axis cs:-0.10800002515316,-1.24252831935883)
--(axis cs:-0.0760000348091125,-1.00236678123474)
--(axis cs:-0.0440000295639038,-0.756528854370117)
--(axis cs:-0.0120000317692757,-0.515527546405792)
--(axis cs:0.0199999660253525,-0.286729007959366)
--(axis cs:0.0519999638199806,-0.0727817714214325)
--(axis cs:0.0839999616146088,0.125859975814819)
--(axis cs:0.115999966859818,0.307051777839661)
--(axis cs:0.147999957203865,0.46556892991066)
--(axis cs:0.179999947547913,0.594264984130859)
--(axis cs:0.211999952793121,0.685999631881714)
--(axis cs:0.24399995803833,0.73531174659729)
--(axis cs:0.275999963283539,0.739395380020142)
--(axis cs:0.307999938726425,0.698349714279175)
--(axis cs:0.339999943971634,0.614844381809235)
--(axis cs:0.371999949216843,0.493481397628784)
--(axis cs:0.403999924659729,0.340139895677567)
--(axis cs:0.435999929904938,0.161466330289841)
--(axis cs:0.467999935150146,-0.0354481935501099)
--(axis cs:0.499999940395355,-0.24323758482933)
--(axis cs:0.531999945640564,-0.454380929470062)
--(axis cs:0.563999950885773,-0.661306738853455)
--(axis cs:0.595999956130981,-0.856574296951294)
--(axis cs:0.627999901771545,-1.03311252593994)
--(axis cs:0.659999907016754,-1.1845269203186)
--(axis cs:0.691999912261963,-1.30537581443787)
--(axis cs:0.723999977111816,-1.3913506269455)
--(axis cs:0.755999982357025,-1.43936347961426)
--(axis cs:0.787999987602234,-1.44754993915558)
--(axis cs:0.819999992847443,-1.41537308692932)
--(axis cs:0.851999998092651,-1.3440535068512)
--(axis cs:0.883999943733215,-1.23734652996063)
--(axis cs:0.915999948978424,-1.10260486602783)
--(axis cs:0.947999954223633,-0.951777815818787)
--(axis cs:0.979999959468842,-0.801335334777832)
--(axis cs:1.01199996471405,-0.669523298740387)
--(axis cs:1.04399991035461,-0.571118831634521)
--(axis cs:1.07599997520447,-0.513702988624573)
--(axis cs:1.10799992084503,-0.498257458209991)
--(axis cs:1.13999998569489,-0.521744549274445)
--(axis cs:1.17199993133545,-0.57893979549408)
--(axis cs:1.2039999961853,-0.663384795188904)
--(axis cs:1.23599994182587,-0.768006265163422)
--(axis cs:1.26800000667572,-0.885636508464813)
--(axis cs:1.29999995231628,-1.00948870182037)
--(axis cs:1.29999995231628,2.06460523605347)
--(axis cs:1.29999995231628,2.06460523605347)
--(axis cs:1.26800000667572,1.96687960624695)
--(axis cs:1.23599994182587,1.83457255363464)
--(axis cs:1.2039999961853,1.66752183437347)
--(axis cs:1.17199993133545,1.46814382076263)
--(axis cs:1.13999998569489,1.24169707298279)
--(axis cs:1.10799992084503,0.996329128742218)
--(axis cs:1.07599997520447,0.742780327796936)
--(axis cs:1.04399991035461,0.493383884429932)
--(axis cs:1.01199996471405,0.259850561618805)
--(axis cs:0.979999959468842,0.049962043762207)
--(axis cs:0.947999954223633,-0.134039491415024)
--(axis cs:0.915999948978424,-0.292916357517242)
--(axis cs:0.883999943733215,-0.426429897546768)
--(axis cs:0.851999998092651,-0.531766712665558)
--(axis cs:0.819999992847443,-0.604436576366425)
--(axis cs:0.787999987602234,-0.639964938163757)
--(axis cs:0.755999982357025,-0.635309278964996)
--(axis cs:0.723999977111816,-0.589667439460754)
--(axis cs:0.691999912261963,-0.504617929458618)
--(axis cs:0.659999907016754,-0.383751332759857)
--(axis cs:0.627999901771545,-0.232107102870941)
--(axis cs:0.595999956130981,-0.0556027591228485)
--(axis cs:0.563999950885773,0.13930931687355)
--(axis cs:0.531999945640564,0.345812737941742)
--(axis cs:0.499999940395355,0.556771397590637)
--(axis cs:0.467999935150146,0.764747262001038)
--(axis cs:0.435999929904938,0.962083578109741)
--(axis cs:0.403999924659729,1.1411120891571)
--(axis cs:0.371999949216843,1.29448628425598)
--(axis cs:0.339999943971634,1.41562008857727)
--(axis cs:0.307999938726425,1.49910879135132)
--(axis cs:0.275999963283539,1.54107928276062)
--(axis cs:0.24399995803833,1.53936886787415)
--(axis cs:0.211999952793121,1.49358296394348)
--(axis cs:0.179999947547913,1.40519857406616)
--(axis cs:0.147999957203865,1.27785742282867)
--(axis cs:0.115999966859818,1.1179713010788)
--(axis cs:0.0839999616146088,0.935551404953003)
--(axis cs:0.0519999638199806,0.744957208633423)
--(axis cs:0.0199999660253525,0.5645672082901)
--(axis cs:-0.0120000317692757,0.413845837116241)
--(axis cs:-0.0440000295639038,0.307981014251709)
--(axis cs:-0.0760000348091125,0.254113435745239)
--(axis cs:-0.10800002515316,0.252056956291199)
--(axis cs:-0.140000015497208,0.296977400779724)
--(axis cs:-0.172000020742416,0.381446242332458)
--(axis cs:-0.204000025987625,0.496665835380554)
--(axis cs:-0.236000016331673,0.633264183998108)
--(axis cs:-0.26800000667572,0.782045543193817)
--(axis cs:-0.300000011920929,0.934599936008453)
--cycle;
\addlegendimage{
      confidence box={circle,fill,inner sep=1pt}{diamond,fill,inner sep=1pt},
      blue}

\addplot [semithick, blue]
table {%
-0.300000011920929 -0.60244745016098
-0.26800000667572 -0.644212424755096
-0.236000016331673 -0.668025732040405
-0.204000025987625 -0.668787717819214
-0.172000020742416 -0.642094612121582
-0.140000015497208 -0.584744930267334
-0.10800002515316 -0.495235681533813
-0.0760000348091125 -0.374126672744751
-0.0440000295639038 -0.224273920059204
-0.0120000317692757 -0.0508408546447754
0.0199999660253525 0.138919115066528
0.0519999638199806 0.336087703704834
0.0839999616146088 0.530705690383911
0.115999966859818 0.712511539459229
0.147999957203865 0.871713161468506
0.179999947547913 0.999731779098511
0.211999952793121 1.0897912979126
0.24399995803833 1.13734030723572
0.275999963283539 1.14023733139038
0.307999938726425 1.09872925281525
0.339999943971634 1.01523220539093
0.371999949216843 0.893983840942383
0.403999924659729 0.740625977516174
0.435999929904938 0.561774969100952
0.467999935150146 0.364649534225464
0.499999940395355 0.156766891479492
0.531999945640564 -0.0542840957641602
0.563999950885773 -0.260998725891113
0.595999956130981 -0.456088542938232
0.627999901771545 -0.632609844207764
0.659999907016754 -0.784139156341553
0.691999912261963 -0.904996871948242
0.723999977111816 -0.990509033203125
0.755999982357025 -1.0373363494873
0.787999987602234 -1.04375743865967
0.819999992847443 -1.0099048614502
0.851999998092651 -0.937910079956055
0.883999943733215 -0.831888198852539
0.915999948978424 -0.697760581970215
0.947999954223633 -0.542908668518066
0.979999959468842 -0.375686645507812
1.01199996471405 -0.204836368560791
1.04399991035461 -0.0388674736022949
1.07599997520447 0.114538669586182
1.10799992084503 0.249035835266113
1.13999998569489 0.359976291656494
1.17199993133545 0.444602012634277
1.2039999961853 0.502068519592285
1.23599994182587 0.533283174037933
1.26800000667572 0.54062157869339
1.29999995231628 0.527558207511902
};
\end{axis}

\end{tikzpicture}}
    \caption{(Left) A \ac{gp} prior with zero mean and \ac{se} covariance function. (Right) The same \ac{gp}, but conditioned on datapoints (black asterisk).
    The blue line is its mean and the blue area is two times its standard deviation ($2\sigma$).}
\end{figure}

\section{Linear Ordinary Differential Equation GPs}\label{sec:LODE_GP}
The class of \acp{gp} is closed under linear operations\footnote{This holds true for almost all relevant cases in control theory. For more details see \citep{harkonen2023gaussian, matsumoto2024images}.}, i.e.\ applying a linear operator $\mathcal{L}$ to a \ac{gp} $g$ as $\mathcal{L}g$ is again a \ac{gp}.
This ensures that realizations of the \ac{gp} $\mathcal{L}g$ lie in the image of the linear operator $\mathcal{L}$, in addition to the \ac{gp} $g$ \citep{jidling2017linearly,langehegermann2018algorithmic}.

We demonstrate the procedure from \citep{besginow2022constraining} for constructing so-called \acp{lodegp}, i.e.\ \acp{gp} that strictly satisfy the underlying system of linear homogenuous \acp{ode}.

\begin{figure}
    \begin{tikzpicture}[scale=0.95]
  \def\tankWidth{2}
  \def\tankHeight{2}
  \def\leftTankPipeDistance{1.5*\tankHeight}
  \def\leftMinPipeDistance{\leftTankPipeDistance-0.5}
  \def\rightTankPipeDistance{1.2+\leftTankPipeDistance}
  \def\rightMinPipeDistance{\rightTankPipeDistance-0.5}
  \def\pipeWidth{0.2*\tankWidth}
  \coordinate (left_tank_x_center) at (0,0);
  \coordinate (middle_tank_x_center) at ($(left_tank_x_center)+(1.2*\tankWidth,0)$);
  \coordinate (right_tank_x_center) at ($(middle_tank_x_center)+(1.2*\tankWidth,0)$);

  \draw ($(left_tank_x_center)+(-0.5*\tankWidth,0)$) arc [
        start angle=0,
        end angle=180,
        x radius=-0.5*\tankWidth,
        y radius=-0.1*\tankHeight
    ]-- ($(left_tank_x_center)+(0.5*\tankWidth,\tankHeight)$)
    arc [
        start angle=0,
        end angle=540,
        x radius=0.5*\tankWidth,
        y radius=-0.1*\tankHeight
    ] --($(left_tank_x_center)+(-0.5*\tankWidth,0)$) ;

  \draw ($(middle_tank_x_center)+(-0.5*\tankWidth,0)$) arc [
        start angle=0,
        end angle=180,
        x radius=-0.5*\tankWidth,
        y radius=-0.1*\tankHeight
    ]-- ($(middle_tank_x_center)+(0.5*\tankWidth,\tankHeight)$)
    arc [
        start angle=0,
        end angle=540,
        x radius=0.5*\tankWidth,
        y radius=-0.1*\tankHeight
    ] --($(middle_tank_x_center)+(-0.5*\tankWidth,0)$);

  \draw ($(right_tank_x_center)+(-0.5*\tankWidth,0)$) arc [
        start angle=0,
        end angle=180,
        x radius=-0.5*\tankWidth,
        y radius=-0.1*\tankHeight
    ]-- ($(right_tank_x_center)+(0.5*\tankWidth,\tankHeight)$)
    arc [
        start angle=0,
        end angle=540,
        x radius=0.5*\tankWidth,
        y radius=-0.1*\tankHeight
    ] --($(right_tank_x_center)+(-0.5*\tankWidth,0)$);

  \draw ($(left_tank_x_center)+(-0.7*\tankWidth, \leftTankPipeDistance)$)
      -- ($(left_tank_x_center)+(-0.5*\pipeWidth, \leftTankPipeDistance)$)
      -- ($(left_tank_x_center)+(-0.5*\pipeWidth, \leftMinPipeDistance)$)
    arc [
        start angle=0,
        end angle=180,
        x radius=-0.5*\pipeWidth,
        y radius=-0.2*\pipeWidth
    ] -- ($(left_tank_x_center)+(0.5*\pipeWidth, \leftTankPipeDistance)$)
      -- ($(middle_tank_x_center)+(-0.5*\pipeWidth, \leftTankPipeDistance)$)
      -- ($(middle_tank_x_center)+(-0.5*\pipeWidth, \leftMinPipeDistance)$)
    arc [
        start angle=0,
        end angle=180,
        x radius=-0.5*\pipeWidth,
        y radius=-0.2*\pipeWidth
    ] -- ($(middle_tank_x_center)+(0.5*\pipeWidth, \pipeWidth+\leftTankPipeDistance)$)
      -- ($(left_tank_x_center)+(-0.7*\tankWidth, \pipeWidth+\leftTankPipeDistance)$);

  \draw ($(right_tank_x_center)+(0.7*\tankWidth, \rightTankPipeDistance)$)
      -- ($(right_tank_x_center)+(0.5*\pipeWidth, \rightTankPipeDistance)$)
      -- ($(right_tank_x_center)+(0.5*\pipeWidth, \rightMinPipeDistance)$)
    arc [
        start angle=0,
        end angle=180,
        x radius=0.5*\pipeWidth,
        y radius=-0.2*\pipeWidth
    ] -- ($(right_tank_x_center)+(-0.5*\pipeWidth, \rightTankPipeDistance)$)
      -- ($(middle_tank_x_center)+(0.5*\pipeWidth, \rightTankPipeDistance)$)
      -- ($(middle_tank_x_center)+(0.5*\pipeWidth, \rightMinPipeDistance)$)
    arc [
        start angle=0,
        end angle=180,
        x radius=0.5*\pipeWidth,
        y radius=-0.2*\pipeWidth
    ] -- ($(middle_tank_x_center)+(-0.5*\pipeWidth, \pipeWidth+\rightTankPipeDistance)$)
      -- ($(right_tank_x_center)+(0.7*\tankWidth, \pipeWidth+\rightTankPipeDistance)$);
  \node (left_tank_text) at ($(left_tank_x_center)+(0, 0.1)$) {$x_1(t)$};
  \node (middle_tank_text) at ($(middle_tank_x_center)+(0, 0.1)$) {$x_2(t)$};
  \node (right_tank_text) at ($(right_tank_x_center)+(0, 0.1)$) {$x_3(t)$};
  \node (left_pipe_text) at ($(left_tank_x_center)+(-0.9*\tankWidth, \leftTankPipeDistance + 0.5*\pipeWidth)$) {$u_1(t)$};
  \node (right_pipe_text) at ($(right_tank_x_center)+(0.9*\tankWidth, \rightTankPipeDistance + 0.5*\pipeWidth)$) {$u_2(t)$};
    \end{tikzpicture}
    \caption{A sketch of the three tank system.}
    \label{fig:three_tank_sketch}
\end{figure}
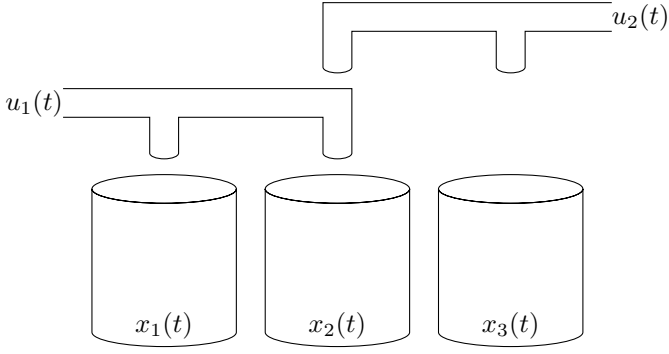
We introduce the \ac{lodegp} approach following the three tank system, illustrated in Figure~\ref{fig:three_tank_sketch}.  
\begin{equation}\label{eq:three_tank_diffeq}
\mathbf{0} 
 = \underbrace{\begin{bmatrix}
-\partial_t & 0 & 0 & 1 & 0 \\
0 & -\partial_t & 0 & 1 & 1 \\
0 & 0 & -\partial_t & 0 & 1
\end{bmatrix}}_A
    \cdot
     \underbrace{\begin{bmatrix}
        x_1(t)\\
        x_2(t)\\
        x_3(t)\\
        u_1(t)\\
        u_2(t)
    \end{bmatrix}}_F
\end{equation}
We can factor $A$ using the \ac{snf} as $UAV = D$, with $D \in \mathbb{R}[\partial_t]^{m_1 \times m_0}$ the one-dimensional representation of the system, $U \in \mathbb{R}[\partial_t]^{m_1 \times m_1}$ and $V \in \mathbb{R}[\partial_t]^{m_0 \times m_0}$ the base change matrices \citep{smith1862systems, newman1997smith}.
All matrices belong to the polynomial ring $\mathbb{R}[\partial_t]$ i.e.\ containing polynomials of $\partial_t$.
For the system in Equation~\eqref{eq:three_tank_diffeq} the application of an algorithm to find the \ac{snf} results in the following $D$, $U$ and $V$:
\begin{align*}
&\underbrace{\begin{bmatrix}
1 & 0 & 0 \\
-1 & 1 & 0 \\
1 & -1 & 1
\end{bmatrix}}_U
\underbrace{\begin{bmatrix}
-\partial_t & 0 & 0 & 1 & 0 \\
0 & -\partial_t & 0 & 1 & 1 \\
0 & 0 & -\partial_t & 0 & 1
\end{bmatrix}}_A
\underbrace{\begin{bmatrix}
0 & 0 & 0 & -1 & 0 \\
0 & 0 & 0 & 0 & -1 \\
0 & 0 & 1 & 1 & -1 \\
1 & 0 & 0 & -\partial_t & 0 \\
0 & 1 & 0 & \partial_t & -\partial_t
\end{bmatrix}}_V
\\=&
\underbrace{\begin{bmatrix}
1 & 0 & 0 & 0 & 0 \\
0 & 1 & 0 & 0 & 0 \\
0 & 0 & -\partial_t & 0 & 0
\end{bmatrix}}_D
\end{align*}
We then construct a prior latent \ac{gp} $\tilde{g}$ using simple construction rules based on the diagonal entries of $D$, which are $0$ or $1$ for controllable systems. 
This construction is based on appropriate basis functions $w$ that result in the corresponding covariance functions, see Appendix D in \citep{besginow2022constraining}.
In the case of the system in Equation~\eqref{eq:three_tank_diffeq}, we construct the latent \ac{gp} $\tilde{g}$: 
\begin{equation}
    \tilde{g} = \mathcal{GP}\left(\begin{pmatrix}0 \\ 0 \\ 0 \\ 0 \\ 0 \end{pmatrix}, 
        \begin{pmatrix}
            0 & 0 & 0 & 0 & 0\\
            0 & 0 & 0 & 0 & 0\\
            0 & 0 & 1 & 0 & 0\\
            0 & 0 & 0 & k_\text{SE} & 0\\
            0 & 0 & 0 & 0 & k_\text{SE}
        \end{pmatrix}\right)
    \end{equation}
which we simplify to $g = \mathcal{GP}\left(\begin{pmatrix}
    0\\0\\0
\end{pmatrix}, \diag\left(\begin{pmatrix}
    1\\
    k_\text{SE}\\
    k_\text{SE}
\end{pmatrix}\right)\right)$ and simplify $V$ similarly by omitting the first two columns. 
By applying the linear base change operator $V$ to this \emph{latent} \ac{gp} $g$ we modify our functional prior $g$ to a different function prior $Vg$ (see Equation~\eqref{eq:LODE_GP_harmonic}) with the property that it \emph{strictly} satisfies the system in Equation~\eqref{eq:three_tank_diffeq}, as detailed in \cite{besginow2022constraining}.
By doing so we guarantee that the \ac{lodegp} $Vg$ spans the nullspace of $A$, which is equivalent to saying that $Vg$ produces only solutions to the original homogenuous system $A\cdot F=\mathbf{0}$.
\begin{equation}\label{eq:LODE_GP_harmonic}
    Vg = \mathcal{GP}\left(\mathbf{0}, V\cdot 
        \diag\left(\begin{pmatrix}
			1\\
			k_\text{SE}\\
			k_\text{SE}
		\end{pmatrix}\right)
    \cdot \hat{V}^T\right)
\end{equation}
where $\hat{V}$ is the operator $V$ applied to the second argument ($t'$) of the \ac{se} covariance function $k_{\text{SE}}$ (cf.\ Equation~\eqref{eq:SE_kernel}).

The \ac{snf} $D$ has the nice property, that its entries are only polynomials with real or complex zeros, or the integers $0$ and $1$, see \cite[Table 1]{besginow2022constraining}. By leveraging this, we obtain the following property about the solution space of a \ac{lodegp} as in \eqref{eq:LODE_GP_harmonic}. 

\begin{table*}[t]
\caption{Primary operators $d$ and their corresponding covariance function $k(t_1, t_2)$. Taken from (\cite{besginow2022constraining}).}
    \centering
\begin{tabular}{cc}
    \hline
    $d$ & $k(t_1, t_2)$\\
    \hline
    $1$ & $0$ \\
    $(\partial_t -a)^j$ & $\left(\sum_{i=0}^{j-1} t_1^it_2^i\right)\cdot \exp(a\cdot (t_1+ t_2))$ \\
    $((\partial_t -a - ib)(\partial_t -a + ib))^j$ & $\left(\sum_{i=0}^{j-1} t_1^it_2^i\right)\cdot \exp(a\cdot (t_1 + t_2))\cdot\cos(b\cdot( t_1-t_2))$ \\
    $0$ & $\exp(-\frac{1}{2}(t_1-t_2)^2)$ \\
    \hline
\end{tabular}
\label{tab:base_cov_construction_paper}
\end{table*}

\begin{thm}
    The set of realizations of a \ac{lodegp} is dense in
    \begin{equation}
        0 \oplus \mathcal{F}^{m_0} \oplus w\cdot\mathbb{R}^q 
    \end{equation}
   with $q$ the number of uncontrollable decoupled systems related to polynomial entries in $D$ and $w$ a vector of basis functions appropriate to describe the non-controllable part for systems of linear homogenuous \acp{ode} with constant coefficients.
\end{thm}
\begin{pf} We consider each component of the solution space $\mathbf{0} \oplus \mathcal{F}^{m_0} \oplus w\cdot\mathbb{R}^q$ individually:

For $\mathbf{0}$: The zero functions trivially belong to this space. The space of realizations of a \ac{gp} with zero mean and zero covariance is equal to the space $\mathbf{0}$ and is therefore also dense.

For $\mathcal{F}^{m_0}$: The denseness of realizations of \acp{gp} with a strong inductive bias for controllable systems of differential equations in \( \mathcal{F}^{m_0}\) has been established in \citep{langehegermann2018algorithmic}.

For $w\cdot\mathbb{R}^q$: 
We choose an appropriate set of basis functions $h$.
\begin{enumerate}
    \item $\mathbb{R}^q$ represents the coefficients of the basis functions defined by the $q$ independent (potentially multiple) zeroes of diagonal entries $d_i \in D$
    \item We assume the distribution $\mathcal{N}(0, I)$ over the parameter space $\mathbb{R}^q$ since $\supp(\mathcal{N}(0, I)) = \mathbb{R}^q$
    \item We apply the basis functions $w$ as $\supp(w_*\mathcal{N}(0, I)) = w\cdot\mathbb{R}^q$
\end{enumerate}
Via this construction we get the prior latent \ac{gp} $\tilde{g}$.
Applying the linear transformation $V$ and $\hat{V}^T$ induces the corresponding isomorphisms between the solution spaces of $A$ and $D$.

Thus, the realizations of the \ac{lodegp} are dense in \( 0 \oplus \mathcal{F}^{m_0} \oplus w\cdot\mathbb{R}^q \).
\qed
\end{pf}

The \ac{lodegp} in Equation~\eqref{eq:LODE_GP_harmonic} can be trained and conditioned on datapoints, same as the regular \ac{gp} in Figure~\ref{fig:GP_example}, which we exploit for our optimal control algorithm.

\section{LODE-GP Optimal Control}\label{sec:method}

This section describes how to apply our \ac{lodegp} for optimal control problems to generate the control input using \ac{gp} conditioning, i.e. control as inference.
We assume a \ac{lodegp} as described in Section~\ref{sec:LODE_GP} for a system of linear \acp{ode}. As we have seen in the previous section, the \ac{lodegp} posterior mean yields smooth functions which satisfy \eqref{eq:mpc:con:ode}. In the following we will show how we enforce constraint \eqref{eq:mpc:con:init} by a specific construction of the conditioned dataset $\mathcal{D}$ and obtain optimality as in \eqref{eq:mpc:obj}.

In order to force a point constraint~\eqref{eq:mpc:con:init} in timestep $i$, we use $(t_i, f_i)$ in the conditioned dataset $\mathcal{D}$ with noise variance $\sigma_n^2(t_i) = 0 \in \mathbb{R}^{m_0}$.  This results in the \ac{lodegp} posterior mean satisfying $\mu^*(t_i) = f_i$ up to numerical precision if it is admissible, i.e. the covariance matrix is invertible with vanishing signal noise $\sigma_n^2$.
Tracking a desired reference trajectory can be given by 
providing further setpoints $(t_i, f_i)$.
In the absence of correlated point constraints which means that the respective second summand in the posterior forumulation \eqref{eq:gaussian_process_posterior_distribution} converges to zero, the posterior (mean and variance) converges to the prior (mean and variance), i.e. the first summand in \eqref{eq:gaussian_process_posterior_distribution} respectively. In \cite{tebbe2024physics}, the authors leverage this fact in order to proof open-loop stability for a Model Predictive Control scheme resulting from optimal control problems as in \eqref{eq:mpc:obj} - \eqref{eq:mpc:con:init}.

For the optimality in \eqref{eq:mpc:obj}, we introduce the concept of \acp{rkhs} \citep{aronszajn1950theory}, which defines a space of functions spanned by the covariance function of our \ac{lodegp}.
\begin{defn}
Let $k$ be a covariance function on $\mathbb{R}$. A Hilbert space $\mathcal{H}_k$ of functions on $\mathbb{R}$ with an inner product $\langle \cdot, \cdot \rangle_{\mathcal{H}_k}$ is called a \ac{rkhs} with reproducing kernel $k$ if 
\begin{enumerate} \label{def:rkhs}
    \item \label{rkhs:bounded} $k(\cdot, t) \in \mathcal{H}_k$ for all $t \in \mathbb{R}$
    \item \label{rkhs:rep_prop}For all $h \in \mathcal{H}_k$ and $t \in \mathbb{R}$ it holds:\begin{center}
        $h(t) = \langle h, k(\cdot, t)\rangle_{\mathcal{H}_k}$
    \end{center} 
\end{enumerate}
with \eqref{rkhs:rep_prop} being the reproducing property.
\end{defn}
Note that each covariance function of a \ac{gp} is also a kernel in the sense an \ac{rkhs}. To this end, we will use both terms equivalently for the remainder of this paper.
 The Moore-Aronszajn Theorem \citep{aronszajn1950theory} yields that there exists a unique \ac{rkhs} for each covariance function (kernel) $k$. 
We can construct the \ac{rkhs} $\mathcal{H}_k$ for a given kernel $k$ from a pre-Hilbert space 
\begin{equation}
    \mathcal{H}_0 = \text{span}\{k(\cdot, t) : t \in \mathbb{R}\}
\end{equation}
with the inner product 
\begin{equation}
    \langle h_1, h_2 \rangle_{\mathcal{H}_0} := \sum_{r=1}^{n} \sum_{s=1}^{m} a_r b_s k(t_r, t'_s)
\end{equation}
for $h_1 = \sum_{r=1}^{n} k(\cdot, t_r)$ and $h_2 = \sum_{s=1}^{m} k(\cdot, t'_s)$. 
The norm $\Vert h \Vert_{\mathcal{H}_0}$ is induced by the inner product and the \ac{rkhs} ${\mathcal{H}_k}$ is defined as the closure of ${\mathcal{H}_0}$ with respect to $\Vert \cdot \Vert_{\mathcal{H}_0}$, i.e.
\begin{align*}
    {\mathcal{H}_k} = \{ h =& \sum_{r=1}^\infty c_r k(\cdot, t_r) : c_r \in \mathbb{R}, t_r \in \mathbb{R}, \\ 
    s.t. \Vert h \Vert_{\mathcal{H}_k}^2 :&= \lim\limits_{n \to \infty} \left\Vert \sum_{r=1}^n c_r k(\cdot, t_r) \right\Vert_{\mathcal{H}_0}^2 \\ &= \sum_{r,s=1}^\infty c_r c_s k(t_r, t_s) < \infty \}.
\end{align*}
The functions of the \ac{rkhs} thereby inherit certain properties of the kernel $k$ such as smoothness. Moreover, a lower \ac{rkhs} norm indicates a smoother function \citep{kanagawa2018gaussian}.
The connection between the mean of a \ac{gp} and its \ac{rkhs} is given via the representer theorem \citep{scholkopf2001generalized}.

\begin{thm}
    For the optimization problem 
    \begin{equation}
        \min\limits_{h \in \mathcal{H}_k} \frac{1}{n} \sum_{i=1}^n (h(t_i) - f_i)^2 + \sigma_n^2 \Vert h \Vert_{\mathcal{H}_k}^2,
    \end{equation}
    with an \ac{rkhs} $\mathcal{H}_k$, the solution is given via
    \begin{equation}
        h^*(t^*) = k(t^*, t)(k(t, t) + \sigma_n^2 I_n)^{-1}f
    \end{equation}
    with $f_i \in \mathbb{R}^{m_0}$.
\end{thm}
This result states, that the problem of finding the best data fit in the (potentially infinite dimensional) \ac{rkhs} using regularization is given in closed form as the weighted sum of kernel evaluations of the \emph{finite} amount of datapoints. For \acp{gp}, the coefficients of this weighted sum yield the formula for the posterior mean in \eqref{eq:gaussian_process_posterior_distribution}. This means that the posterior mean is in some sense optimal in the \ac{rkhs} of functions given the covariance function $k$. 
If the regularization factor $\sigma_n^2$ is zero, the posterior mean minimizes the \ac{rkhs} norm and fits the datapoints $f$ exactly.

\begin{thm}\label{thm:rkhs}
    The posterior mean function as stated in Equation~\eqref{eq:gaussian_process_posterior_distribution} with noise variance $\sigma_n^2 = 0$ yields an optimal control, minimizing the \ac{rkhs} norm given by the covariance function $k$ of the LODE-GP with the constraints \eqref{eq:mpc:con:ode} and \eqref{eq:mpc:con:init}, i.e. the optimization problem
    \begin{equation}
        \min\limits_{h \in \mathcal{H}_k} \Vert h \Vert_{\mathcal{H}_k} \text{ s.t. } h(t_i) = f_i, h \in \text{sol}_\mathcal{F}(A)
    \end{equation}
\end{thm}
\begin{pf}
    Since $\text{sol}_\mathcal{F}(A)$ is closed in $C^\infty$, it is also closed in $\mathcal{H}_k$.
	Then the result follows directly from \cite[Theorem 3.5]{kanagawa2018gaussian} and the fact that a \ac{lodegp} satisfies $h \in \text{sol}_\mathcal{F}(A)$. \qed
\end{pf}
This result corresponds to the objective~\eqref{eq:mpc:obj} and provides the optimality of our optimal control problem. It completes the translation of the optimization problem described in chapter~\ref{sec:optimalcontrol} to a problem solved by inference of a \ac{lodegp} with a specific training dataset. Note that the \ac{rkhs} is generally dependent on the hyperparameters of the kernel. In the following two sections we will investigate how the solutions of the optimization problem \eqref{eq:mpc:obj} - \eqref{eq:mpc:con:init} differ dependent on the hyperparameters which are lengthscale and signal variance for the \ac{se} kernel in \eqref{eq:SE_kernel}.

\begin{figure*}[t]
    \centering
     \includegraphics[width=0.335\linewidth]{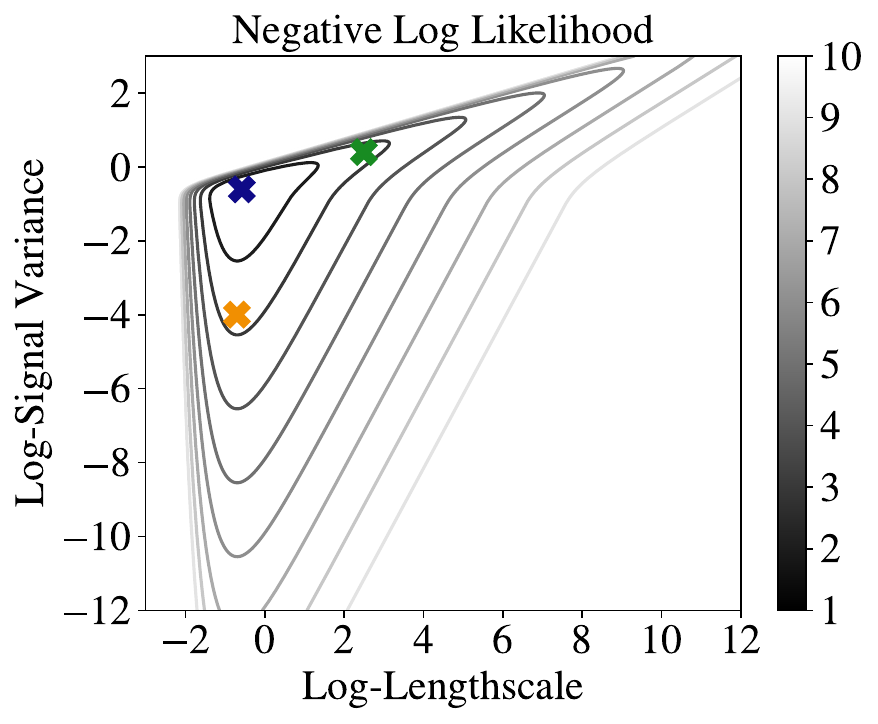}
    \scalebox{0.5}{\input{plots/minimal_opt.tex}}
    \scalebox{0.5}{\input{plots/minimal_lhigh_s_high2.tex}}
    \scalebox{0.5}{\input{plots/minimal_llow_s_low2.tex}}
    
    \caption{Minimal system: In the left plot we see a contour plot of the negative log likelihood with three colored crosses representing the hyperparamter combinations $(\ell^2, \sigma_f^2)$ on log scale. We pick the optimal combination (blue) and two further hyperparameter combinations (orange and green) to demonstrate the influence of them on the optimal control problem. We plot the mean and two times the standard deviation as shaded area.}
    \label{fig:minimal}
\end{figure*}
\section{Example I}\label{sec:evaluation}

First, we evaluate a minimal system in two variants. The first variant requires a minimal number of states, i.e. one state:

\begin{equation*}
    A_1 = \begin{pmatrix} 1-\partial_t
\end{pmatrix},
D_1 = \begin{pmatrix}
     1
\end{pmatrix},
U_1 = \begin{pmatrix}
     1 
\end{pmatrix},
V_1 = \begin{pmatrix}
    1 
\end{pmatrix}.
\end{equation*}

By adding a control on the derivative, the system, its base change matrices and \ac{snf} become:

\begin{equation*}
    A_2 = \begin{pmatrix} 1 & \partial_t
\end{pmatrix},
D_2 = \begin{pmatrix}
     1 & 0\\
     0 & 0
\end{pmatrix},
U_2 = \begin{pmatrix}
     1 
\end{pmatrix},
V_2 = \begin{pmatrix}
     0 & 1\\
     -1 & \partial_t
\end{pmatrix}.
\end{equation*}

This allows us to perform optimal control, since we have one degree of freedom. 

If we want to steer the system from a certain state, e.g. $x_0 = (1,0)^T$ to the origin $(0,0)^T$, we obtain the optimization problem 
\begin{align*}\label{ex:oc:minimal}
    &\min \Vert F \Vert_\mathcal{F} \\
    \text{s.t. } &F \in \text{sol}_\mathcal{F}(A_2) \\
    &F(0) = (1,0)^T \\
    &F(1) = (0,0)^T
\end{align*}

solved by the corresponding \ac{lodegp} conditioned on $\mathcal{D} = \{(0, (1, 0)^T), (1, (0, 0)^T)\}$.
Since we have an \ac{se} kernel on the state channel, we have two essential hyperparameters, the lengthscale $\ell^2$ and the signal variance $\sigma_f^2$. Recall that we do not vary the noise variance in contrast to standard \ac{gp} regression.  The lengthscale indicates the bandwidth of setpoints. In our example, a small lengthscale results in a fast decay, see the rightmost and leftmost graphs in Figure \ref{fig:minimal} of the state towards the prior. On the other hand, a larger lengthscale would result in a slow decay, see the center graph in Figure \ref{fig:minimal}. The optimal hyperparameter combination based on the optimized \ac{mll} is $(\ell^2_*, \sigma_{f*}^2) \approx (\exp(-0.5774), \exp(-0.6097))$, see the elevation plot in Figure \ref{fig:minimal}. For the solution of the optimal control problem, we may choose the mean as argumented in the last section. For further tasks as determining safety probabilities \citep{tebbe2024efficiently}, we may use the whole distribution given via mean and posterior covariance.

\section{Example II: Three Tank System}

\begin{figure*}[h]
    \centering
    \scalebox{0.46}{\input{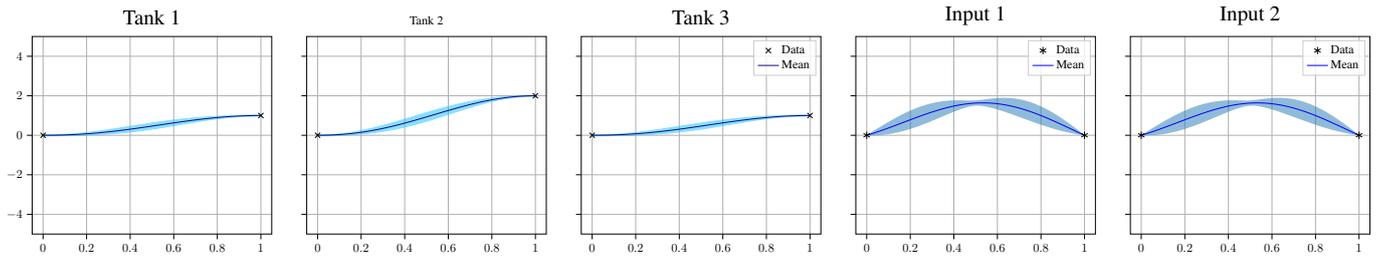}}
    \caption{Plot of the three tank system for the given optimal control problem with optimized hyperparameters based on the \ac{mll}. The signal variance of the uncontrollable part is close to zero since the optimal control problem is admissible.}
    \label{fig:three_tank}
\end{figure*}

Furthermore, we have the three tank system containing two kernels and one uncontrollable part, arising from the states being linearly dependent, see Section~\ref{sec:LODE_GP}. Due to the algorithm of the \ac{snf}, the uncontrollable part is generated as the state of the third tank. We model the uncontrollable part as a normal distribution with signal variance $\sigma_u^2$. The value of this is trained to nearly 0 if the optimal control problem is admissible for the given dataset $\mathcal{D}$. Note that the matrices $U$ and $V$ are not unique. We observe this behavior also for an optional $V_3$ with 
\begin{equation}
    V_3 = \begin{bmatrix}
0 & 0 & 1 & 1 & 0 \\
0 & 0 & -2 & 1 & 1 \\
0 & 0 & 1 & 0 & 1 \\
1 & 0 & 0 & \partial_t & 0 \\
0 & 1 & 0 & 0 & \partial_t
\end{bmatrix}
\end{equation}
and respective matrix $U$. 
For a similar optimal control task as for the minimal system, we obtain
\begin{align*}\label{ex:oc:threetank}
    &\min \Vert F \Vert_\mathcal{F} \\
    \text{s.t. } &F \in \text{sol}_\mathcal{F}(A) \\
    &F(0) = (0,0,0,0,0)^T \\
    &F(1) = (1,2,1,0,0)^T
\end{align*}

as the problem of steering the regulated system to a steady and admissible state.
The solution for optimized hyperparameters is presented in Figure \ref{fig:three_tank}. Note that it is debatable to train all hyperparameters based on the setpoints. The uncontrollable part yields an inner state of the three tanks which cannot be broken, so this state should be learned from historic data rather than from future manually chosen setpoints. 
\section{Conclusion}
\label{sec:conclusion}
We have discussed an approach for solving optimal control problems in smooth function spaces using computer algebra combined with \acp{gp}. We apply a strong inductive bias, forcing the \ac{gp} to only produce functions following the given differential equations. Providing the pointwise constraints of the optimal control problem as conditioning data to the \ac{gp}, we obtain a solution for the optimal control problem directly via \ac{gp} inference. Our method defined by a \ac{lodegp} posterior mean yields optimal functions in the norm of the corresponding \ac{rkhs}, therefore reducing an optimal control problem to an inference problem.

\bibliography{sample}
\end{document}